\documentclass[b5paper,reqno]{amsart}
\usepackage{amsfonts,lingmacros,tree-dvips, amsmath, amssymb, graphicx, mathrsfs}
\usepackage[hidelinks]{hyperref}

\usepackage[all]{xy}

\newdir{ >}{!/6pt/@{ }*:(1,0.2)@_{>}*:(1,-0.2)@^{>}}

\theoremstyle{plain}

\newtheorem{theorem}{Theorem}[section]

\theoremstyle{definition}
\newtheorem{definition}[theorem]{Definition}

\newtheorem{remark}[theorem]{Remark}
\newtheorem{counter example}[theorem]{Counter Example}

\numberwithin{equation}{section}

\DeclareMathAlphabet{\mathscr}{OT1}{pzc}{m}{it} 
\begin{document}
\Large{
		\title{A NOTE ON SETS AVOIDING RATIONAL DISTANCES IN CATEGORY BASES}
		
		\author[S. Basu]{Sanjib Basu}
		\address{\large{Department of Mathematics,Bethune College,181 Bidhan Sarani}}
		\email{\large{sanjibbasu08@gmail.com}}
		
		\author[A.C.Pramanik]{Abhit Chandra Pramanik}
	\address{\large{Department of Pure Mathematics, University of Calcutta, 35, Ballygunge Circular Road, Kolkata 700019, West Bengal, India}}
	\email{\large{abhit.pramanik@gmail.com}}

		\thanks{The second author thanks the CSIR, New Delhi – 110001, India, for financial support}
	\begin{abstract}
   Michalski gave a short and elegant proof of a theorem of A. Kumar which states that for each set $A\subseteq\mathbb{R}$, there exists a set $B\subseteq A$ which is full in $A$ and such that no distance between points in $B$ is a rational number. He also proved a similar theorem for sets in $\mathbb{R}^2$. In this paper, we generalize these results in some special types of category bases.
	\end{abstract}
\thanks{}
	\maketitle

\section{INTRODUCTION}
In Euclidean $n$-dimensional space $\mathbb{R}^n$, a set which is such that the distance between any two of its points is always an irrational number is called an ARD set, or, ``a set avoiding rational distance". Komjath [$4$] asked the following question : Do every set in $\mathbb{R}^n$ contains an ARD set such that both the sets have the same Lebesgue outer measure ? In [$5$], A. Kumar gave a positive answer to this question for linear sets. Here is a statement Kumar's theorem as given by Michalski in his paper [$6$].
\begin{theorem}
	Let $A\subseteq\mathbb{R}$ be a set of positive outer measure. Then there exists an ARD set $B\subseteq A$ full in $A$.\\
	(Here by the phrase ``$B$ is full in $A$" we mean that env($A$)=env($B$), where for any set $X$, env($X$) is a $G_\delta$ set containing $X$ such that the inner measure of env($X)-X$ is zero)
\end{theorem}
Kumar's proof depends on a theorem of Githik and Shelah [$2$] and use some results on forcing with $\sigma$-ideals which are quite complicated. However, Michalski gave a short and elegant proof of Theorem $1.1$. As the relation of being at rational distance between points in $\mathbb{R}^2$ is not an equivalence relation, this approach cannot be replicated for planar sets. But using an altogether different method, Michalski proved that 
\begin{theorem}
	If $A$ be a measurable subset of $\mathbb{R}^2$ of positive measure, then there exists an ARD set $B\subseteq A$ full in $A$.
\end{theorem}

Also, under the assumption that non($\mathcal{N}$) = cof($\mathcal{N}$), where $\mathcal{N}$ denotes the family of Lebesgue null sets, he established a relatively stronger form of the above theorem.
\begin{theorem}
	Let $A\subseteq \mathbb{R}^2$ be a measurable set of positive measure. Then there exists a partial bijection $f$ the graph of which is an ARD full subset of $A$.
\end{theorem}
It may be noted that on account of the similarities between the $\sigma$-algebras of measurable sets and sets having Baire property and between the $\sigma$-ideals of measure zero sets and sets of first category, it is possible to give category analogues of the above two theorems where the statement ``$B$ is full in $A$" may be replaced by ``$B$ is of second category in $A$" which means that $X\cap B$ is of second category whenever $X\cap A$ is so.\\

The approach of Michalski is our main motivation to look for some generalization of the above theorems in some special types of category bases, where each of the above theorems along with their category analogous can be unified under a common framework.

\section{PRELIMINARIES AND RESULTS}
The concept of a category base is a generalization of both measure and topology. Its main objective is to present measure and Baire category (topology) and also some other aspects of point set classification within a common framework. It was introduced by J. C. Morgan II [$7$] in the seventies of the last century and since then has been developed through a series of papers [$1$], [$8$], [$9$], [$10$], [$11$] etc. 

\begin{definition}
A pair $(X,\mathcal{C}$) where $X$ is a non-empty set and $\mathcal{C}$ is a family of subsets of $X$ is called a category base if the non-empty members of $\mathcal{C}$ called regions satisfy the following set of axioms :
\begin{enumerate}
\item Every point of $X$ belongs to some region; i.e., $X = \cup$ $\mathcal{C}$.
\item Let $A$ be a region and $\mathcal{D}$ be any non-empty family of disjont regions having cardinality less than the cardinality of $\mathcal{C}$.\\
i) If $A \cap ( \cup \mathcal{D}$) contains a region, then there exists a region $D\in\mathcal{D}$ such that $A\cap D$ contains a region. \\
ii)  If $A\cap(\cup \mathcal{D})$ contains no region, then there exists a region $B\subseteq A$ which is disjoint from every region in $\mathcal{D}$.
\end{enumerate}
\end{definition}

Several examples of category bases are given in [$7$]

\begin{definition}
In a category base ($X,\mathcal{C}$), a set $A$ is called `singular' if every region contains a subregion which is disjoint from the set. A set which can be expressed as countable union of singular sets is called `meager'. Otherwise, it is called `abundant'. A set whose complement is meager is called `co-meager'. A set is `Baire' (or, a Baire set) if every region contains a subregion in which the set or its complement is meager. 
\end{definition}
\begin{theorem}
	(The Fundamental Theorem) In a category base ($X,\mathcal{C}$), every abundant set is abundant everywhere in some region. This means that for every abundant set $A$, there exists a region $C$ in every subregion $D$ of which $A$ is abundant.
\end{theorem}
To the above list of Definitions, we further add some special types of category bases.

\begin{definition}
A category base is called `point-meager' if every singleton set in it is meager, `separable' if there is a countable subfamily of regions such that every region is abundant everywhere in at least one region in the subfamily, `perfect' if $X=\mathbb{R}^n$ and for every region $A$ and for every point $x\in A$, there is a descending sequence \{$A_n\}_{n=1}^{\infty}$ of regions such that $x\in A_n$ and diam($A_n)\leq {1/n}, n=1,2,..$ and a `translation base' if $X=\mathbb{R}$ and 
\begin{enumerate}
	\item $\mathcal{C}$ is translation invariant.
	\item If $A$ is any region and $D$ is a countable everywhere dense set, then $\bigcup\limits_{r\in D}A(r)$ is abundant everywhere, where $A(r)=\{x+r: x\in A\}$
\end{enumerate}
\end{definition}
Following Grzegorek and Labuda [$3$], we may say that 
\begin{definition}
	In a category base ($X,\mathcal{C}$), a set $F$ is a full subset of $E$ if $F\subseteq E$ and for every Baire set $B$, $B\cap F$ is abundant whenever $B\cap E$ is so. If $E$ is a Baire set, this is equivalent to saying that $E-F$ cannot contain any abundant Baire set.
\end{definition} 
Definition $2.5$ provides a common generalization for two analogous concepts of full subset in measure and category within the framework of category bases. Our first two theorems are generalizations of Theorem $1.1$ in perfect bases. In what follows, we use symbol $\mathbb{Q}$ for the set of rationals and $c$ for the cardinality of the continuum. The same symbol $c$ is also used for the smallest ordinal representing it.

\begin{theorem}
	Let $(\mathbb{R},\mathcal{C})$ be any perfect base. Then every abundant Baire set $E$ contains an ARD set $F$ which is a full subset of $E$ in $(\mathbb{R},\mathcal{C})$.
\begin{proof}
Let $\mathcal{F}=\{P: P$ is a perfect set and $P\subseteq E\}$. By (Th 8, II, Ch 5, [7]), $\mathcal{F}\neq \phi$. Again as each perfect set contains $c$ number of disjoint perfect sets (Th $20$, I, Ch $5$, [$7$]), so $\mathcal{F}$ contains at least $c$ number of elements in it. But the cardinality of the family of all perfect sets in $\mathbb{R}$ is also $c$. Therefore card$(\mathcal{F})=c$.\\
We write $\mathcal{F}=\{P_\alpha: \alpha<c\}$ and construct a subset $F$ of $E$ by defining an injective family $\{x_\alpha:\alpha<c\}$ of points in $E$ in the following manner: suppose for an ordinal $\alpha$, we have already defined the partial family $\{x_\beta:\beta<\alpha\}$. Now consider the set $Z_\alpha=\bigcup\{x_\beta+\mathbb{Q}:\beta<\alpha\}$ and choose $P_\alpha\in\mathcal{F}$. Since card$(Z_\alpha)\leq$ card$(\alpha).\omega<c$ whereas card$(P_\alpha)=c$, so $P_\alpha-Z_\alpha\neq\phi$ and we choose $x_\alpha\in P_\alpha-Z_\alpha$ and proceed accordingly. We finally set $F=\{x_\alpha:\alpha<c\}$.\\
The set $F$ is obviously an ARD set. If possible, let there be an abundant Baire set $B$ such that $B\cap E$ is abundant but $B\cap F$ is meager in $(\mathbb{R},\mathcal{C})$. Then $B\cap (E-F)$ is an abundant Baire set and so contains a perfect set $T$. Certainly, $T\in \mathcal{F}$ and so accordingly to our construction $T$ must meet $F$ non-vacously. But, this is impossible. Hence $F$ is a full subset of $E$ in $(\mathbb{R},\mathcal{C})$.\\ 
This proves the theorem.	
	\end{proof}
\end{theorem}
However, an alternative proof of the above theorem can be given using Vitali-Bernstein selector (i.e. a selector in the quotient group $\mathbb{R}/\mathbb{Q}$ which is both Vitali as well as a Bernstein set). Let $V$ be a Vitali-Bernstein selector. Regarding construction of such a selector, we refer to Th $2.2$ [$6$]. Since $E=\bigcup\limits_{r\in\mathbb{Q}}((V+r)\cap E)$, so there exists $r_0\in\mathbb{Q}$ such that $(V+r_0)\cap E$ is abundant. We set $F=(V+r_0)\cap E$. The set $F$ being a subset of a Vitali set is obviously an ARD set. Moreover, it is a full subset of $E$ , for otherwise, it would be possible to fit an abundant Baire set in $E$ which is disjoint with $F$. But this is impossible, because by (Th $8$, Ch $5$, [$7$]) every abundant Baire set in any perfect base contains a perfect set and $V$ being Bernstein should intersect this perfect set.
\begin{remark}
	We are not certain whether in the above theorem we can remove the restriction of Baireness from the set $E$. But this may be done as the next result shows if our perfect base is constitued of all perfect sets in $\mathbb{R}$ [$7$]. As in the standard terminology [$7$], here also we write Marczewski meager, Marczewski abundant, Marczewski Baire for meager, abundant and Baire sets in the category base of all perfect sets.
\end{remark}
\begin{theorem}
	Let $(\mathbb{R},\mathcal{P})$ be the category base of all perfect sets in $\mathbb{R}$ and $E$ is Marczewski abundant. Then there exists an ARD set $F$ which is a full subset of $E$ in $(\mathbb{R},\mathcal{P})$.
	\begin{proof}
		Let $\mathcal{E}=\{P:P\in\mathcal{P}$ such that card($E\cap P)=c$ \}. By the Fundamental theorem there exists $Q\in\mathcal{P}$ such that $E$ is Marczewski abundant everywhere in $Q$. But this implies that for every perfect set $D$ contained in $Q$, card($E\cap D)=c$, for otherwise, there exists a perfect set $T\subseteq Q$ such that card($E\cap T)<c$. But then by (Th $20$, I, Ch $5$, [$7$]) there exists a perfect set $T^\prime\subseteq T\subseteq Q$ which is disjoint with $E$ - a contradiction. Hence $\mathcal{E}$ contains at least $c$ number of elements in it. Again the cardinality of the family of  all perfect sets in $\mathbb{R}$ is also $c$. Therefore card($\mathcal{E})=c$.\\
		We write $\mathcal{E}=\{P_\alpha:\alpha<c\}$ and construct a subset $F$ of $E$ by defining an injective family \{$x_\alpha:\alpha<c\}$ of points in $E$ as follows : suppose for an ordinal $\alpha$, we have already defined the partial family \{$x_\beta:\beta<\alpha$\}. Now consider the set $Z_\alpha=\bigcup\{x_\beta+\mathbb{Q}:\beta<\alpha$\} and choose $P_\alpha\in\mathcal{E}$. Since card($Z_\alpha)\leq$card($\alpha$)$.\omega<c$ and card($E\cap P_\alpha)=c$, so $(E\cap P_\alpha)-Z_\alpha\neq\phi$. We further choose $x_\alpha\in(E\cap P_\alpha)-Z_\alpha$ and proceed accordingly. Finally, we set $F=\{x_\alpha:\alpha<c\}$.\\
		The set $F$ is obviously an ARD set. If possible, let there be a Marczewski abundant Baire set $B$ such that $B\cap E$ is Marczewski abundant but $B\cap F$ is Marczewski meager. By the Fundamental theorem and (Th $8$, Ch $5$, [$7$]), there exists a perfect set $Q$ contained in $B$ such that such that $B\cap E$ is Marczewski abundant everywhere in $Q$. Consequently, from the above construction it follows that every perfect subset of $Q$ contains at least one member of $F$. But $B\cap F$ being Marczewski meager in $(\mathbb{R},\mathcal{P})$, there exists at least one perfect subsets of $Q$ which is disjoint with $F$ - a contradiction. Hence $F$ is a full subset of $E$.\\
			This proves the theorem.
		\end{proof}
\end{theorem}
\begin{remark}
	In addition to being a perfect base, if the category base $(\mathbb{R},\mathcal{C})$ in Theorem $2.6$ is also a translation base, then the set $F$ cannot be Baire. This is established once we show that no abundant subset of $F$ is Baire. If possible, let there be an abundant Baire set $H$ contained in $F$. Choose $r^\prime\in\mathbb{Q}$ and fix it. Since $\mathbb{Q}-r^\prime$ is dense in $\mathbb{R}$, $\bigcup\limits_{r\in\mathbb{Q}-\{r^\prime\}}(H+r)$ is comeager (Th $3$, Sec I, Ch $6$, [$7$]). Consequently, the set $H+r^\prime$ which is disjoint with $\bigcup\limits_{r\in\mathbb{Q}-\{r^\prime\}}(H+r)$ is meager implying that $H$ is meager - a contradiction.\\
	But $(\mathbb{R},\mathcal{P})$ is not a translation base for no translation base can be equivalent to the category base of all perfect sets in $\mathbb{R}$ (Note following Th $3$, I, Ch $6$, [$7$]). So the set $F$ in Theorem $2.8$ cannot be proved as non-Baire in the same way as above. But it can be extended to a non-Baire set. Let $\hat{\mathcal{E}}=\{P:P\in\mathcal{P}$ such that card($P\cap E)<c\}$. For each $P\in\hat{\mathcal{E}}$, we choose an element from $P-E$ and thereby extend $F$ to a set $B$ such that $B\cap E=F$. Certainly, $B$ is a Bernstein set which is non-Baire in any perfect base. 
\end{remark}
Let $(X,\mathcal{C})$ and $(Y,\mathcal{D})$ be two category base and $(X\times Y, \mathcal{C}\times \mathcal{D})$ be their product (assume that the product make sense). We denote the product $\mathcal{C}\times \mathcal{D}$ by the symbol $\mathcal{M}$. If $(X,\mathcal{C})$ and $(Y,\mathcal{D})$ are separable, then so is $(X\times Y,\mathcal{M})$ (Th $5.5$, [$10$]). Again if $(X,\mathcal{C})$ and $(Y,\mathcal{D})$ are Baire bases, then so is $(X\times Y,\mathcal{M})$ (Th $5.3$, [$10$]). Hence $(X\times Y,\mathcal{M})$ satisfies countable chain condition (CCC) (Th $2$, I, Ch $3$, [$7$]) which implies that every Baire set can be expressed as the union of a $\mathcal{M}_{\sigma\delta}$-set and a meager set (Th $3$, III, Ch $1$, [$7$]). The following theorem which is a formulation of Theorem $1.3$ in category bases uses this notion of separability in product bases. As Theorem $1.3$ is stronger than Theorem $1.2$, generalizing the later also serves the purpose of generalizing the former.
\begin{theorem}
	Let $(\mathbb{R},\mathcal{C})$ and $(\mathbb{R},\mathcal{D})$ be two separable Baire bases where card$(\mathcal{C})= card(\mathcal{D})=k$ such that $\aleph_0$ is not cofinal with respect to $k$. Moreover, suppose any subset of $\mathbb{R}$ having cardinality less than $k$ is both $\mathcal{C}$ as well as $\mathcal{D}$-meager. Then every abundant Baire set $E$ in $(\mathbb{R}^2,\mathcal{M})$ contains a partial bijection $f$ the graph of which is an ARD full subset of $E$.
	\begin{proof}
	Here card$(\mathcal{M})=k$ and therefore card$(\mathcal{M}_{\sigma\delta})=k^{\aleph_0}$	which is same as $k$ because $\aleph_0$ is not cofinal with $k$. We may assume that the class of sets from $\mathcal{M}_{\sigma\delta}$ whose intersection with $E$ is abundant has the cardinality $k$, for otherwise, we may repeat a set in $\mathcal{M}_{\sigma\delta}$ infinitely often. Consider an enumeration \{$M_\alpha:M_\alpha\in\mathcal{M}_{\sigma\delta}$ and $M_\alpha\cap E$ is abundant in $(\mathbb{R}^2,\mathcal{M})$\} of this class.\\
	 Let $N_\alpha=M_\alpha\cap E$ $(\alpha<k)$ (Here we use the same symbol $k$ also to denote the smallest ordinal representing the cardinal $k$). Then there exists $A_\alpha$ which is abundant in $(\mathbb{R},\mathcal{C})$ such that the $x$-sections $(N_\alpha)_x$ of $N_\alpha$ are abundant in  $(\mathbb{R},\mathcal{D})$ for every $x\in A_\alpha$, for otherwise, $N_\alpha$ would be meager in $(\mathbb{R}^2,\mathcal{M})$ (Th $5.2$, [$10$]) - a contradiction.\\
	 We now construct the desired function $f$ by transfinite induction using procedure similar to Michalski in the present categorical settings. Suppose for any $\alpha<k$, we have already defined the family \{$(x_\xi,y_\xi):\xi<\alpha$\} which is an ARD subset of $E$ and a partial bijection satisfying $(x_\xi,y_\xi)\in N_\xi$. Since every set of cardinality less than $k$ is both $\mathcal{C}$ as well as $\mathcal{D}$-meager, we can choose $x_\alpha\in A_\alpha-\{x_\xi:\xi<\alpha\}$ and $y_\alpha\in(N_\alpha-\bigcup\limits_{\xi<\alpha}\bigcup\limits_{q\in\mathbb{Q}}S((x_\xi,y_\xi), q))_{x_\alpha}$ (where $S((x_\xi,y_\xi), q)$ represents the set of all points in $\mathbb{R}^2$ which are at a distance $q$ from $(x_\xi,y_\xi)$) so that we can further extend the partial bijection whose graph is \{$(x_\xi,y_\xi):\xi<\alpha$\} by adjoining the point $(x_\alpha,y_\alpha)\in N_\alpha$ to it. This completes the construction of $f$.\\
	 Obviously $f$ is an ARD subset of $E$. Moreover, it is a full subset of $E$ in $(\mathbb{R}^2,\mathcal{M})$ because every abundant Baire set in $(\mathbb{R}^2,\mathcal{M})$ which meets $E$ abundantly contains a set from the class \{$N_\alpha:\alpha<k$\}.\\
	 Hence the theorem.
		\end{proof}
\end{theorem}
\begin{remark}
	Since by hypothesis, every singleton set in $\mathbb{R}$ is $\mathcal{D}$-meager, $f$ must be non-Baire in $(\mathbb{R}^2,\mathcal{M})$ for otherwise $f$ would be meager (Th $5.2$, [$10$]) and consequently $E$ would also be meager - a contradiction.
\end{remark}
\bibliographystyle{plain}

	\end{document}